\documentclass[lettersize,journal]{IEEEtran}
\usepackage{amsmath,amsfonts}
\usepackage{algorithmic}
\usepackage{array}
\usepackage[caption=false,font=normalsize,labelfont=sf,textfont=sf]{subfig}
\usepackage{textcomp}
\usepackage{stfloats}
\usepackage{url}
\usepackage{verbatim}
\usepackage{graphicx}

\usepackage{balance}

\usepackage{mathtools} 

\usepackage{amssymb} 

\newcommand{\QED}{\hfill\ensuremath{\blacksquare}} 

\newcounter{constnum}
\newcommand{\const}[1]{\refstepcounter{constnum}\label{#1}}

\newtheorem{theorem}{Theorem}
\newtheorem{lemma}{Lemma}
\newtheorem{proposition}{Proposition}

\newtheorem{remark}{Remark}

\newtheorem{corollary}{Corollary}

\usepackage[style=ieee, backend=bibtex]{biblatex}
\usepackage{customCommands}

\usepackage[caption=false,font=normalsize]{subfig}

\usepackage{enumerate}
\usepackage{xcolor}
\begin{document}
\title{Non-Asymptotic Analysis of Classical Spectrum Estimators for $L$-mixing Time-series Data with Estimated Means}
\author{Yuping Zheng, Andrew Lamperski
\thanks{This work was supported by NSF ECCS 2412435 and NIH 1R01NS147767-0.
}
\thanks{Y. Zheng and A. Lamperski are with the department of Electrical and Computer Engineering, University  of Minnesota, Minneapolis, MN 55455, USA {\tt\small zhen0348@umn.edu, alampers@umn.edu}}}


\maketitle

\begin{abstract}
Spectral estimation is an important tool in time series analysis, with applications including economics, astronomy, and climatology. The asymptotic theory for non-parametric estimation is well-known but the development of non-asymptotic theory is still ongoing. 
Our recent work obtained the first non-asymptotic error bounds on the Bartlett and Welch methods with restrictive assumptions. In this work, we derive non-asymptotic error bounds for both Bartlett and Welch estimators for $L$-mixing time-series data with unknown means, and the results cover the special case with known zero means. 
The class of $L$-mixing processes contains common models in time series analysis, including autoregressive processes and measurements of geometrically ergodic Markov chains. Our new error bounds are of $O(\frac{1}{\sqrt{k}})$, where $k$ is the number of data segments used in the algorithm. Such bounds are the tightest among the existing work on non-asymptotic analysis of classical spectrum estimators with or without zero-mean assumptions.


\end{abstract}


\section{Introduction}

Spectral estimation is widely used in signal processing and control systems. The central objective of spectrum estimation is to estimate the power spectral density of a signal, given a finite data record. Spectrum estimators can be primarily classified into parametric and non-parametric methods \cite{stoica2005spectral}.
Parametric methods include ARMA and state space models. Fitting parametric model requires knowledge or  assumptions about the process,  such as autoregressive order.
Non-parametric methods include periodograms, the Blackman-Tukey method, the Bartlett method, and the Welch method. Non-parametric approaches can be more robust when little is known about the target system, since only the definition of power spectral density is used \cite{stoica2005spectral}.

The asymptotic analysis of spectral estimation has reached a substantial level of maturity (see, e.g. \cite{brockwell1991time, brillinger2001time, stoica2005spectral, liu2010asymptotics, zhang2021convergence}). 
Non-asymptotic theory aims to remove the assumption of asymptotically many data samples, but is less developed. For parametric estimators,
non-asymptotic analysis has been widely studied. In particular, the non-asymptotic results for autoregressive models have been established since the early 2000s \cite{goldenshluger2001nonasymptotic}. In addition, in dynamic system identification, non-asymptotic error bounds are presented in \cite{hardt2018gradient,sarkar2019finite,oymak2019non,tsiamis2019finite,lee2020non}. In contrast, non-asymptotic analysis for non-parametric estimators is less developed. Existing non-asymptotic analysis of non-parametric spectral estimators includes the work \cite{fiecas2019spectral} on smoothed periodograms, variations on Blackman-Tukey estimators in \cite{zhang2021convergence,veedu2021topology}, and \cite{doddi2022efficient} on Wiener filters. Our prior work on non-asymptotic spectral estimation includes \cite{lamperski2023nonasymptotic}, which gives a framework for deriving error bounds for a family of estimators, under rather restrictive assumptions, and \cite{zheng2025non}, which gives error bounds for the Bartlett and Welch algorithms under the assumption of zero-mean $L$-mixing data. 

Similar to our previous paper, \cite{zheng2025non}, we focus on $L$-mixing time-series data.
%
$L$-mixing processes cover many models in time series analysis including most practical cases satisfying the assumptions in \cite{lamperski2023nonasymptotic}, autoregressive processes, and measurements of uniformly geometrically ergodic Markov chains \cite{gerencser2002new}.
The class of $L$-mixing processes quantifies the decay of dependencies of stochastic processes over time and was first introduced in \cite{gerencser1989class}.
Other related work uses the theory of $L$-mixing processes to study stochastic optimization algorithms with time-correlated data streams \cite{chau2021stochastic,zheng2022constrained}. 
Time dependencies have also been described by other mixing conditions (see, e.g. \cite{doukhan1994mixing}), but the application of these conditions to spectral  analysis is beyond the scope of this paper. 

Our contribution is to obtain non-asymptotic error bounds for Bartlett and Welch estimators for $L$-mixing time-series data with unknown means for both batch and online algorithms. 
This work extends the results in \cite{zheng2025non} from zero-mean assumption to the general case that mean is unknown. The main results in this work include the case of zero-mean dataset and the error bounds are tighter than the prior bounds obtained in work \cite{lamperski2023nonasymptotic,zheng2025non}.  
In particular, these are the tightest error bounds for the Bartlett and Welch algorithms, whether the mean is known or unknown. 
As is standard \cite{stoica2005spectral}, the batch method centers the data using the sample mean.
The online method iteratively uses segments of data to estimate the mean and the power spectral density, and is suitable for both streaming data and large data sets.
This work represents a stepping stone for the non-asymptotic analysis of non-parametric spectral estimation for other less ideal time-dependent data. 

 This work unifies and substantially extends the preliminary frameworks introduced in our previous conference paper \cite{zheng2025non} and arXiv preprint \cite{zheng2025non_b}. In particular, we synthesize the distinct cases involving both known and unknown means into a cohesive set of unified results. Furthermore, due to conference page constraints, our prior works could only offer brief sketches of the proofs for the extensions of $L$-mixing properties. To bridge this gap, this article provides the complete, rigorous mathematical proofs omitted in those short-form publications, establishing a definitive theoretical foundation.

Despite making progress in developing non-asymptotic error bounds for spectrum estimators, the works \cite{fiecas2019spectral, lamperski2023nonasymptotic, zheng2025non} are all limited to zero-mean data.  
Although the zero-mean assumption is widely adopted since the data can be demeaned \cite{stoica2005spectral, brillinger1981time}, subtracting an incorrect mean (which would arise when using the sample mean) will lead to estimation errors. Therefore, conducting the non-asymptotic analysis under the assumption that the mean of the dataset is unknown is of practical importance. The work \cite{dudek2023spectral} uses a subsampling method to avoid data demeaning and handle non-stationary data, but only covers the asymptotic analysis for almost periodically correlated time series. While it would be interesting to generalize our theory to non-stationary data, this is beyond the scope of the current paper.

The paper is organized as follows. In Section \ref{sec:Setup}, we describe the problem and algorithms, give background on the class of $L$-mixing processes. Section \ref{sec:convergence_analysis} presents the main results on spectral estimation error analysis. Section~\ref{sec:extension_Lmixing} shows the extensions of classical $L$-mixing results and the detailed proofs. In Sections \ref{sec:proofs_supporting} and \ref{sec:proofs_main}, the proofs of supporting lemmas and the main results are shown, respectively. In Section  \ref{sec:sim}, we verify our theory through simulations of a finite-state Markov chain. Conclusions are given in Section~\ref{sec:conclusion}.
\section{Problem Setup} \label{sec:Setup}

\subsection{Notation}
Random variables are denoted in bold. If $\bx$ is a random variable, then 
$\bbE[\bx]$ is its expected value. 
$\bbR$, $\bbC$, $\bbN$, and $\bbZ$ denote the set of real numbers, complex numbers, nonnegative integers, and integers. 
For a vector, $x$, $\|x\|_2$ denotes the Euclidean norm.
For a matrix, $A$, $\|A\|_F$ denotes the Frobenius norm.
$A^\star$ denotes the conjugate transpose of $A$. We denote the floor function of $x \in \bbR$ by $\floor{x} \coloneqq \max\{n \in \bbZ \vert n \le x\}$.

Let $\cY$ be a finite-dimensional vector space with inner product $\langle \cdot,\cdot\rangle$ and corresponding norm $\|\cdot \|$. For a random variable, $\by\in\cY$, and $q\ge 1$ let $\|\by\|_{L_q}=\left(\bbE[\|\by\|^q]\right)^{1/q}$, which is the corresponding $L_q$ norm. Function $f_t \in L_p[0,T]$ implies that $(\int_0^T |f_t|^p dt)^{\frac{1}{p}} < \infty$.

%
%

\subsection{Problem and Algorithm}

Given a stationary discrete-time stochastic process $\by[k] \in \bbR^n$, the mean, autocovariance sequence, and power spectral density are
\begin{align*}
\mu &= \bbE[\by[i]] \\
R[k] &= \bbE[(\by[i+k] - \mu) (\by[i] - \mu)^\top] \\
\Phi(s) &= \sum_{k = -\infty }^{\infty} e^{-j 2 \pi s k} R[k]
\end{align*}
where $s \in [-1/2, 1/2]$. 

Assuming the true mean $\mu$ is unknown, we want to analyze two non-parametric spectrum estimators, the Bartlett and the Welch estimators. Both methods rely on decomposing the data into segments of $M$ data points and working with Fourier transform estimates. For streamlined notation, the sample mean and Fourier transform estimates for the $i$th data segment are denoted by:
\begin{align*}
\bar{\by}_i &= \frac{1}{M} \sum_{k=0}^{M-1} \by[iK+k]\\
              \hat{\by}_i(s) &= \sum_{k=0}^{M-1} w_k(s) \by[iK +k].
\end{align*}
The Bartlett method has $K=M$ and $w_k(s)=\frac{1}{\sqrt{M}}e^{-j2\pi ks}$, while the Welch method has $K \le M$ and uses $w_k(s)=\frac{v_k}{\|v\|_2}e^{-j2\pi ks}$ for some window vector $0\ne v\in\bbR^M$. Note that when $K<M$, as in the Welch method, the data segments overlap. In both the Bartlett and Welch estimators, $w(s)$ is a Euclidean unit vector for all $s$: $\sum_{k=0}^{M-1}|w_k(s)|^2=1$. In general, we will assume that $K\le M$. (If $K>M$, then the estimators will lose data.)

The running average of the data segment sample means is given by:
$$\hat{\bmu}_k = \frac{1}{k} \sum_{i=0}^{k-1} \bar{\by}_i.$$

Note that
$\bbE[\hat\by_i(s)]=\sum_{i=0}^{M-1}h(s)\mu$, where
$$h(s) = \sum_{k=0}^{M-1} w_k(s).$$


Set $
\Delta \hat{\by}_k(s) = \hat{\by}_k - h(s) \hat{\bmu}_k.
$
 The proposed algorithms are the following:
\begin{itemize}
\item \textit{Batch Algorithm:}
\begin{equation}  \label{eq:alg_batch}
\begin{aligned}
\hat \bPhi_k(s)&=\frac{1}{k}\sum_{i=0}^{k-1}  \left(\hat \by_i(s) - h(s) \hat{\bmu}_k \right) \left(\hat \by_i(s) - h(s) \hat{\bmu}_k \right)^\star,
\end{aligned}
\end{equation}
\item \textit{Online Algorithm:} 

\begin{equation} \label{eq:alg_online}
\begin{aligned}
\hat{\bmu}_{k+1} &= \hat{\bmu}_k + \alpha_k \left(\bar{\by}_k - \hat{\bmu}_k \right) \\
\hat{\bPhi}_{k+1}(s) & = \hat{\bPhi}_{k}(s) + \alpha_k \left( \Delta \hat \by_k(s) \Delta \hat \by_k(s)^\star - \hat{\bPhi}_k(s)\right),
\end{aligned}
\end{equation}
\end{itemize}
where $\alpha_k=\frac{1}{k+1}$.

The batch algorithm, in particular, reduces to the standard presentations of the Bartlett and Welch algorithms from \cite{stoica2005spectral}, after subtracting the mean from the data. The online algorithm provides real-time responsiveness and memory efficiency. Also, when the mean is known and zero, the batch and online algorithms reduce to the batch-form estimators and the iterative algorithm respectively in \cite{zheng2025non}.

Our convergence results cover the general cases of dataset with both unknown mean and known mean and the special case of zero-mean dataset.
When the mean is known, it is equivalent to say that the perfect mean estimate exists. 
Then we can always subtract such known mean value from the data to obtain a zero-mean dataset. This is the main motivation that many existing works simply assume the data is zero-mean \cite{zheng2025non,lamperski2023nonasymptotic}.  
\subsection{Background on the Class of $L$-mixing Processes} \label{subsec:L-mixing_intro}

In this work, we assume that $\by[k]$ is an $L$-mixing data sequence. In this subsection, we present some background on the class of $L$-mixing processes. We start with the classical definitions in continuous
time and describe how they change for discrete time.

 Let $\cF=(\cF_t)_{t\ge 0}$ be an increasing family of $\sigma$-algebras.
 Let $\cF^+=(\cF_t^+)_{t\ge 0}$ be a decreasing family of $\sigma$-algebras such that for all $t\ge 0$, $\cF_t$ and $\cF_t^+$ are independent,
$\cF_t^+=\cF_0^+$ for all $t\le 0$,
and $\cF_t^+=\sigma\left\{\bigcup_{\epsilon >0} \cF_{t+\epsilon} \right\}$. 
A continuous-time stochastic process $\by_t\in\cY$ is called \emph{$L$-mixing} with respect to $(\cF,\cF^+)$ if 
\begin{itemize}
\item $\by_t$ is measurable with respect to $\cF_t$ for all $t\ge 0$
\item $M_q(\by):=\sup_{t\ge 0}\|\by_t\|_{L_q}<\infty$ for all $q\ge 1$
\item $\Gamma_q(\by):=\int_0^{\infty}\gamma_q(\tau,\by)d\tau<\infty$ for all $q\ge 1$, where $\gamma_q(\tau,\by)=\sup_{t\ge \tau}\|\by_t-\bbE[\by_t|\cF_{t-\tau}^+]\|_{L_q}$. 
\end{itemize}
The number, $\Gamma_q(\by)$ characterizes the speed at which dependencies decay over time.

Now, we introduce the definition of discrete-time $L$-mixing process with some mild modification from the continuous-time definition. Let $\cF = (\cF_k)_{k \ge 0}$ and $ \cF^+ =(\cF_k^+)_{k \ge 0}$ be monotone increasing and decreasing sequences of $\sigma$-algebras, respectively, such that $\cF_k$ and $\cF_k^+$ are independent for all $k \ge 0$. 
%
%
Then the discrete-time process $\by_k$ is $L$-mixing with respect to $(\cF, \cF^+)$ if the continuous-time process defined by $\by_t = \by_{\floor{t}}$ is $L$-mixing with respect to the continuous-time $\sigma$-algebras defined by $\cF_t = \cF_{\floor{t}}$ and $\cF_t^+ = \cF_{\floor{t}}^+$.

The moment bounds remains the same as in continuous time, but the discrete-time counterpart of $\Gamma_q(\by)$ is defined by $\Gamma_{d,q}(\by) = \sum_{k=0}^\infty \gamma_q(k,\by)$. 
Note that 
$\gamma_{q}(\tau,\by)=\max\{\gamma_q(\floor{\tau},\by),\gamma_q(\floor{\tau}+1,\by)\}$, which implies that
$
\Gamma_{d,q}(\by)\le \Gamma_q(\by)\le 2\Gamma_{d,q}(\by).
$
In particular, $\Gamma_{d,q}(\by)$ is finite  if and only if $\Gamma_q(\by)$ is finite.


\begin{remark}
$L$-mixing processes can be defined in both continuous time and discrete time. For the results presented in this work, we examine discrete-time $L$-mixing processes. 
Corollary~\ref{thm:L_mixing_sum} below is a discrete-time, multi-dimensional variation of Theorem 1.1 in \cite{gerencser1989class}, for scalar continuous-time processes.
\end{remark} 

\subsection{Comparison with Assumptions from \cite{lamperski2023nonasymptotic}}
\label{ss:comparison}

In this paper, we assume that $\by[k]$ is an $L$-mixing sequence. In contrast, the work in \cite{lamperski2023nonasymptotic} requires that at least one of following two assumptions hold:
\begin{enumerate}[{A}1)]
\item
  \label{a:gaussian}
  $\by[k]$ is Gaussian
\item
  \label{a:subgaussian}
  There is an impulse  response sequence $h[k]\in\bbR^{n\times m}$ such that
  $\by[k] = \sum_{\ell=-\infty}^{\infty} h[k-\ell]\bzeta[\ell]$, where $\bzeta[k]=\begin{bmatrix}\bzeta_1[k] & \cdots & \bzeta_m[k]\end{bmatrix}^\top$ such that for $i=1,\ldots,m$ and for $k\in \bbZ$, $\bzeta_i[k]$ are independent $\sigma$-sub-Gaussian random variables.
\end{enumerate}


The class of $L$-mixing processes contains a wide variety of processes that cannot be modeled by the assumptions from \cite{lamperski2023nonasymptotic}, including measurements of geometrically ergodic Markov chains, \cite{gerencser2002new}, which can be used to model a wide variety of stable nonlinear stochastic systems. 
\footnote{The work in \cite{gerencser2002new} only proves the $L$-mixing property for measurements of uniformly geometrically ergodic Markov chains satisfying a $1$-step Doeblin condition. Based on our ongoing work, we hypothesize that measurements of irreducible, $V$-uniformly ergodic Markov chains are also $L$-mixing. The class of irreducible, $V$-uniformly ergodic Markov chains is broad, and covers nonlinear stochastic dynamical systems which satisfy a stochastic Lyapunov stability condition. See~\cite{meyn2012markov,douc2018markov}.} 
In particular, measurements of an ergodic Markov chain on a finite state space are $L$-mixing, but do not fit the assumptions from \cite{lamperski2023nonasymptotic}.  On the other hand, in many cases, the processes satisfying the assumptions from \cite{lamperski2023nonasymptotic} are also $L$-mixing. 
%
%
%

Furthermore, the class of $L$-mixing processes remains $L$-mixing under a variety of operations. In particular, 
if $\bzeta[k]$ is $L$-mixing, passing $\bzeta[k]$ through a stable, causal linear filter results in another $L$-mixing sequence. 
(Specific conditions on the filter are discussed in the result below.)
 If $f$ is Lipschitz, then $f(\bzeta[k])$ is $L$-mixing. The product of two $L$-mixing sequences is also $L$-mixing. As a result, rather complex processes can be shown to be $L$-mixing. 

The next result shows under suitable hypotheses on the filter $h$, data satisfying assumption A\ref{a:subgaussian} are also $L$-mixing. It is proved in Subsection V-D in \cite{zheng2025non}.

\begin{proposition}
  \label{prop:comparison}
  {\it
    If $\by[k]$ satisfies A\ref{a:subgaussian} and $h$ is causal with $\sum_{\ell=0}^{\infty}\|h[\ell]\|_2(\ell+1)<\infty$, then $\by[k]$ is $L$-mixing with
    \begin{align*}
      M_q(\by)&\le 8m\sigma \sqrt{q}\sum_{\ell=0}^{\infty} \|h[\ell]\|_{2}\\
\Gamma_{d,q}(\by)&\le 8m\sigma \sqrt{q}\sum_{\ell=0}^{\infty} \|h[\ell]\|_{2}(\ell+1).
    \end{align*}
}
\end{proposition}

\begin{remark}
In the case that $\by[k]$ is Gaussian, under the conditions that $\Phi(s)$ admits a spectral factorization (see \cite{wiener1957prediction}), $\by[k]$ also satisfies A\ref{a:subgaussian} with a causal filter $h$, $m=n$, and $\bzeta[k]$ Gaussian. In particular, when $\by[k]$ is generated by a stable linear Gaussian state space system, $h$ can be computed from the Kalman filter, and the hypotheses of Proposition~\ref{prop:comparison} hold. However, finding more general conditions to ensure that $h$ satisfies the requirements of Proposition~\ref{prop:comparison} is out of the scope of this paper.
\end{remark}

\subsection{$L$-mixing Properties of Transformed Data}

For the sake of notation simplicity, we set $\tilde{\by}(s) = \hat{\by}(s) - h(s) \mu$ throughout the paper. Since the error bounds do not depend on $s$, we drop $s$ in the corresponding proofs for compact notation.

The following lemma describes how $L$-mixing properties of the original data sequence induce $L$-mixing properties  
on the vectors and matrices used in the spectral estimation algorithms. 
\const{c_M_y}
\const{c_Gamma_y}
\begin{lemma} \label{lem:L-mixing_constant_collection}
\it{
If $\by[k]$ is L-mixing with respect to $(\cF, \cF^+)$, then for all $s \in \bbR$, $\tilde{\by}(s)$ and $\tilde{\by}(s) \tilde{\by}^\star(s)$ are L-mixing with respect to $(\cG, \cG^+)$ where $\cG_i = \cF_{iK+M-1}$ for all $i \in \bbN$. Furthermore, for all $q \ge 1$, the following bounds hold:
\begin{align*}
&M_{2q}(\tilde{\by}(s)) \le  4 \sqrt{(2q-1) M_{2q}(\by) \Gamma_{d,2q}(\by)} \coloneqq c_{\ref{c_M_y},q} \\
&\Gamma_{d,2q}(\tilde{\by}(s)) \le \\
& \hspace{18pt} 2 (\floor{\frac{M-1}{K}} + 1)  c_{\ref{c_M_y},q} +  (\floor{\frac{M-1}{K}} + 1) \Gamma_{d,2q}(\by)
 \coloneqq c_{\ref{c_Gamma_y},q}\\
&M_{2q}(\tilde{\by}(s) \tilde{\by}(s)^\star) \le c_{\ref{c_M_y},2q}^2 \\
&\Gamma_{d,2q}(\tilde{\by}(s) \tilde{\by}(s)^\star) \le 6 c_{\ref{c_M_y},2q} c_{\ref{c_Gamma_y},2q}.
\end{align*}
In particular, if $\mu$ is known and $\mu=0$, then we obtain tighter bounds:
\begin{align*}
&M_{2q}(\tilde{\by}(s)) \le  2 \sqrt{2(2q-1) M_{2q}(\by) \Gamma_{d,2q}(\by)} \coloneqq \tilde{c}_{\ref{c_M_y},q} \\
&\Gamma_{d,2q}(\tilde{\by}(s)) \le \\
& \hspace{18pt} 2 (\floor{\frac{M-1}{K}} + 1)  \tilde{c}_{\ref{c_M_y},q} +  (\floor{\frac{M-1}{K}} + 1) \Gamma_{d,2q}(\by)
 \coloneqq \tilde{c}_{\ref{c_Gamma_y},q}\\
&M_{2q}(\tilde{\by}(s) \tilde{\by}(s)^\star) \le \tilde{c}_{\ref{c_M_y},2q}^2 \\
&\Gamma_{d,2q}(\tilde{\by}(s) \tilde{\by}(s)^\star) \le 6 \tilde{c}_{\ref{c_M_y},2q} \tilde{c}_{\ref{c_Gamma_y},2q}.
\end{align*}
}
\end{lemma}

\section{Convergence Analysis} \label{sec:convergence_analysis}

Let $\bar{\Phi}(s) = \bbE[(\hat{\by}_i(s) - h(s) \mu ) (\hat{\by}_i(s) - h(s) \mu)^\star]$ where $\mu = \bbE[\by[i]] = \bbE[\bar{\by}_i] = \bbE[\hat{\bmu}_k]  $ due to the stationary assumption. 

In this section, we present a supporting result followed by the main results. 

\subsection{Convergence of Mean Estimates}

The following lemma quantifies the deviation between $\hat{\bmu}_k$ and $\mu$ in expectation, which helps to obtain the explicit error bounds in the main results. Note that $\hat{\bmu}$ is an unbiased estimator.

\begin{lemma} \label{lem:mean_estimate_error}
\it{
Let $\by[k]$ be an $L$-mixing sequence. Assume $\hat{\bmu}_0 = 0 $, $\alpha_j= \frac{1}{j+1}$ for all integers $j\in [0,k-1]$, then for all integers $k \ge 1$ and all $q\ge 1$: 
$$
\|\hat{\bmu}_k - \mu \|_{L_{2q}} \le  c_{q} \frac{ 1 }{\sqrt{M}} \frac{1}{\sqrt{k}} 
$$
where $c_{q} \hspace{-2pt}=\hspace{-2pt}c_{\ref{c_M_y},q} \left( \floor{\frac{M-1}{K}} + 1 \right)$ and $c_{\ref{c_M_y},q}$ is defined in Lemma~\ref{lem:L-mixing_constant_collection}.
}
\end{lemma}

\subsection{Main Results}
The following theorems show the concentration of the estimates around their expected value for both batch and online algorithms.

\begin{theorem} \label{theorem_main_batch}
  {\it
    Let $\by[k]$ be an $L$-mixing sequence. For the batch algorithm \eqref{eq:alg_batch}
    with total number of iterations, $ k \ge 1$, and all $q\ge 1$:
\begin{align*}
\MoveEqLeft[0]
\left\|\hat{\bPhi}_k(s) - \bar{\Phi}(s)] \right\|_{L_{q}}  \hspace{ -5pt} \\
&\le \left(  4 \sqrt{6 (2q-1) c_{\ref{c_M_y},2q}^3 c_{\ref{c_Gamma_y},2q}} + 2 c_{\ref{c_M_y},2q} c_{2q}   + c_{2q}^2  \right) \frac{1}{\sqrt{k}}.
\end{align*}
where $c_{\ref{c_M_y},2q}$, $c_{\ref{c_Gamma_y},2q}$, and $c_{2q}$ are defined in Lemma~\ref{lem:L-mixing_constant_collection} and Lemma~\ref{lem:mean_estimate_error}.
}
\end{theorem}

\begin{theorem} \label{theorem_main_online}
  {\it
    Let $\by[k]$ be an $L$-mixing sequence. For the online algorithm \eqref{eq:alg_online}, assume that $\hat{\bmu}_0=0$, $\alpha_j = \frac{1}{j+1}$, $\forall j \in \bbN $ and $j \in [0, k-1]$, 
    then for all integers $ k \ge 2$ and all $q\ge 1$:
\begin{align*}
\MoveEqLeft[0]
\left\|\hat{\bPhi}_k(s) - \bar{\Phi}(s)] \right\|_{L_{q}} \\
 &  \le b_{q} \frac{1}{\sqrt{k}} +  (M_{4q}(\by)^2  + 2 c_{\ref{c_M_y},2q}  M_{4q}(\by) ) \frac{M}{k}.   
\end{align*}
where 
$
b_q =  4 \sqrt{6 (2q-1) c_{\ref{c_M_y},2q}^3 c_{\ref{c_Gamma_y},2q}}  + 6c_{\ref{c_M_y},2q} c_{2q} + 2  c_{2q}^2 
$
and $c_{\ref{c_M_y},2q}$, $c_{\ref{c_Gamma_y},2q}$, and $c_{2q}$ are defined in Lemma~\ref{lem:L-mixing_constant_collection} and Lemma~\ref{lem:mean_estimate_error}.
}
\end{theorem}


The following Corollary presents a tighter bound when the mean of the dataset is known:
\begin{corollary} \label{cor:known_mean}
  {\it
Let $\by[k]$ be an $L$-mixing sequence with known mean and $\alpha_j = \frac{1}{j+1}$, $\forall j \in \bbN$ and $j \in [0,k-1]$.  For both the batch and online algorithms, for all integers $k \ge 1$ and all $q \ge 1$:
\begin{align*}
\MoveEqLeft[0]
\left\|\hat{\bPhi}_k(s) - \bar{\Phi}(s)] \right\|_{L_{q}} \le  4 \sqrt{6 (2q-1) c_{\ref{c_M_y},2q}^3 c_{\ref{c_Gamma_y},2q}} \frac{1}{\sqrt{k}}.
\end{align*}
In particular, when $\mu=0$, we obtain a tighter bound:
\begin{align*}
\MoveEqLeft[0]
\left\|\hat{\bPhi}_k(s) - \bar{\Phi}(s)] \right\|_{L_{q}} \le  4 \sqrt{6 (2q-1) \tilde{c}_{\ref{c_M_y},2q}^3 \tilde{c}_{\ref{c_Gamma_y},2q}} \frac{1}{\sqrt{k}}.
\end{align*}
where $c_{\ref{c_M_y},2q}$, $c_{\ref{c_Gamma_y},2q}$, $\tilde{c}_{\ref{c_M_y},2q}$, $\tilde{c}_{\ref{c_Gamma_y},2q}$ are defined in Lemma~\ref{lem:L-mixing_constant_collection}. 
}
\end{corollary}

If the factors from Theorem~\ref{theorem_main_batch}, Theorem~\ref{theorem_main_online}, and Corollary~\ref{cor:known_mean} grow polynomially in $q$, then error bounds with high probability can be implied. The Markov chain we use in Section \ref{sec:sim} and the example in Proposition 1 of \cite{zheng2025non} both satisfy the polynomial growth assumption and thus have the bound in Theorem~\ref{thm:highProb} below. Furthermore, to complete the error bound analysis, the bounds on bias are also given in Proposition~\ref{prop:bias}.  The proofs of Theorem~\ref{thm:highProb} and Proposition~\ref{prop:bias} in this work have minimal change compared with Theorem 2 and Proposition 2 in \cite{zheng2025non}, and thus are omitted.

The bounds from Theorem~\ref{theorem_main_batch},  Theorem~\ref{theorem_main_online}, and Corollary~\ref{cor:known_mean} share a unified form, $f_k g_q= \left(a_1 \frac{1}{\sqrt{k}} + a_2 \frac{M}{k} \right) g_q$. Here, we can choose proper constants $a_1$ and $a_2$ and a  function of $q$, $g_q$, in the three scenarios respectively. The next result shows that if $g_q$ is a polynomial, then high-probability bounds can be obtained. 
\begin{theorem}
  \label{thm:highProb}
  {\it
  If error bounds from Theorem~\ref{theorem_main_batch},  Theorem~\ref{theorem_main_online}, and Corollary~\ref{cor:known_mean} can be bounded above by $f_k q^r$ for all $q \ge 1$ and $r >0$, then for all $\nu\in (0,1)$ and all $k \ge 2$:
    \begin{multline*}
      \bbP \left(
        \|\hat{\bPhi}_k(s)-\bar{\Phi}(s)\|_F> f_k e^r\max\left\{1,\frac{(\ln \nu^{-1})^{r}}{r^r}\right\}
      \right)
      \le \nu.
    \end{multline*}
  }
\end{theorem}

\begin{remark}
For Bartlett and Welch methods, typically, we have $\frac{K}{M} \le 1$ and the total number of data is $N = (k-1)K+ M$, which imply $k \le \frac{N}{K}$. In practice, a commonly used value is $\frac{M}{K} = 2$ for the Welch estimator. Therefore, when the mean is unknown, Theorem~\ref{thm:highProb} gives a bound of $O(\sqrt{\frac{K}{N}})$ for the batch algorithm and a bound of $O(\sqrt{\frac{K}{N}})$ when $k \gg M$ for the online algorithm. When the mean is known, Theorem~\ref{thm:highProb} gives a bound of $O(\sqrt{\frac{K}{N}})$ for both batch and online algorithms. These results are tighter than those in \cite{lamperski2023nonasymptotic, zheng2025non}. Specifically, the bound is of $O(\sqrt{\frac{\log \hat{N}}{N/\hat{N}}})$ in \cite{lamperski2023nonasymptotic} and the bound is of $O\left(\frac{\log_2(\log_2(N/K))}{\sqrt{N/K}}\right)$ in \cite{zheng2025non}. Note $\hat{N} <N$ is a tunable parameter which specify the trade-off between bias and variance \cite{lamperski2023nonasymptotic}.
\end{remark}
\begin{proposition} \footnote{
      The bias bounds here correct some mistakes from \cite{zheng2025non}. Namely, we fixed the typo in the factor from $M_q(\by)$ to $M_2(\by)$ on all of the bias  bounds, and the dependence on the weighting function for the Welch algorithm was corrected from $\sum_{i=|k|}^{M-1}\frac{v_{i-|k|}v_i}{\|v\|_2^2}$ to  $\left(1-\sum_{i=|k|}^{M-1}\frac{v_{i-|k|}v_i}{\|v\|_2^2} \right)$.} \label{prop:bias}
  {\it
  If $\by[k]$ is $L$-mixing then:
    \begin{itemize}
  \item The bias of the Bartlett estimator is bounded by
  \begin{multline*}
    \|\Phi(s)-\bar{\Phi}(s)\|_2\le 2 M_2(\by )\sum_{|k|\ge M}\gamma_2(|k|,\by)+ \\
    \frac{ 2 M_2(\by )}{M}\sum_{|k|<M} |k| \gamma_2(|k|,\by).
  \end{multline*}
  \item The bias of the Welch estimator is bounded by
  \begin{multline*}
    \|\Phi(s)-\bar{\Phi}(s)\|_2\le  2 M_2(\by) \sum_{|k|\ge M}\gamma_2(|k|,\by)+\\
    2 M_2(\by )\sum_{|k|<M}\gamma_2(|k|,\by)\left(1-\sum_{i=|k|}^{M-1}\frac{v_{i-|k|}v_i}{\|v\|_2^2} \right).
  \end{multline*}
  \end{itemize}
  }
\end{proposition}

\begin{remark}
If we further assume that the mean of $\by[k]$ is known and zero, the bounds above are of factor 2 tighter. 
\end{remark}

  \begin{remark}
    Combined the bounds on the deviation of the spectral estimate to its expectation with the bounds on the bias in Proposition~\ref{prop:bias} above, we can see the trade-off between bias and variance. In particular, bounds in Theorem~\ref{theorem_main_online} increases monotonically with $M$, while the bound from Proposition~\ref{prop:bias} decreases monotonically with $M$. 
  \end{remark}

\section{Variations on Classical $L$-mixing Results} \label{sec:extension_Lmixing}

Extensions of classical $L$-mixing results to the case of stochastic processes with values in an arbitrary finite-dimensional inner product space are presented here. The key result, Corollary~\ref{thm:L_mixing_sum}, is crucial to prove Lemma~\ref{lem:L-mixing_constant_collection}, Lemma~\ref{lem:mean_estimate_error}, Theorem~\ref{theorem_main_batch}, Theorem~\ref{theorem_main_online}, and Corollary~\ref{cor:known_mean}.

The following generalizes Lemma 2.3 of \cite{gerencser1989class}.

\begin{lemma}
  {\it
    \label{lem:innerProd}
  If $\by_t\in\cY$ is a zero-mean $L$-mixing process with respect to $(\cF,\cF^+)$ and $\bz\in\cY$ is $\cF_s$-measurable  with $s\le t$, then for any $p\ge 1$ and $q\ge 1$ with $\frac{1}{p}+\frac{1}{q}$, we have
  $$
  \left|\bbE\left[\langle \by_t,\bz\rangle\right] \right|\le 2\gamma_p(t-s,\by)\|\bz\|_{L_q}.
  $$
  }
\end{lemma}
{\it Proof}

  We follow the same steps as the proof  in \cite{gerencser1989class}, but replace multiplication with
 an inner product. The details are given below:

 Let $\by_{t,s}^+ =\bbE\left[\by_t \vert \cF_s^+ \right]$. 
 Decompose as following:
 \begin{align*}
    \bbE\left[\langle \by_t, \bz \rangle\right] &= \bbE\left[ \langle \by_{t,s}^+ + \left( \by_t -  \by_{t,s}^+\right), \bz \rangle  \right] \\
    & = \bbE\left[ \langle \by_{t,s}^+ , \bz \rangle  \right] + \bbE\left[ \langle  \by_t -  \by_{t,s}^+, \bz \rangle  \right] \\
    & = \langle \bbE\left[ \by_{t,s}^+ \right], \bbE\left[ \bz\right]   \rangle + \bbE\left[ \langle   \by_t -  \by_{t,s}^+ , \bz \rangle  \right]
 \end{align*}
where the third equality uses the fact that $\by_{t,s}^+$ is independent of $\cF_s$. Additionally, since $\bbE[\by_t] = 0$, we have
$$
\bbE\left[\by_{t,s}^+ \right] = - \bbE\left[ \by_t -  \by_{t,s}^+ \right]. 
$$
Therefore, taking the absolute value and applying the Cauchy–Schwarz inequality and Jensen's inequality, we obtain 
\begin{align*}
\left|  \bbE\left[\langle \by_t, \bz \rangle\right] \right| &\le \left\| \bbE\left[ \by_t -  \by_{t,s}^+ \right]\right\| \left\|  \bbE\left[ \bz\right] \right\| + \bbE\left[ \left| \langle  \by_t -  \by_{t,s}^+ , \bz \rangle  \right| \right] \\
& \le\| \by_t -  \by_{t,s}^+ \|_{L_p} \|\bz\|_{L_q}   +  \|\by_t -  \by_{t,s}^+\|_{L_p} \|\bz\|_{L_q}
\end{align*}
where for the second inequality the first term uses the monotonicity of $L_p$ norm and the second term uses H\"{o}lder's inequality. Therefore, using the definitions of $\gamma_p(t-s,\by)$ completes the proof.
\QED

The following supporting lemma is used to prove Lemma~\ref{lem:integral} and it is the same as Lemma 2.5 in \cite{gerencser1989class}. The proof is also presented in \cite{gerencser1989class}, so we omit the proof here.
\begin{lemma} \label{lem:support_innerProd}
\it{ 
  Let $g_t, k_t$ be nonnegative measurable functions in $[0,T]$ such that $g_t$ is bounded and for all $0 \le t \le T$ and $m \ge 1$,
$$
g_t \le \int_0^t g_s^{1-\frac{1}{m}} k_s ds.
$$
Then, the following holds 
$$
g_t \le \left( \frac{1}{m} \int_0^t k_s ds \right)^m.           
$$
}
\end{lemma}

The following is the generalization of Theorem 1.1 in \cite{gerencser1989class}.

\begin{lemma}
  \label{lem:integral}
  {\it
    Let $\bu_t\in\cY$ be a zero-mean $L$-mixing process and let $f_t\in \bbC$ and $f_t \in L_2[0,T]$. For all $m\ge 1$ and all $T\ge 0$:
    \begin{align*}
      \left\|\int_0^Tf_t \bu_t dt\right\|_{L_{2m}} \hspace{-18pt} \le 
      2\left((2m-1)M_{2m}(\bu)\Gamma_{2m}(\bu)\hspace{-3pt}\int_0^T \hspace{-3pt} |f_t|^2 dt\right)^{\frac{1}{2}}\hspace{-5pt}.
    \end{align*}    
  }
\end{lemma}
{\it Proof}  


Let $\bx_t = \int_0^t f_s \bu_s ds$ and $g_t = \bbE\left[ \|\bx_t\|^{2m}\right]$.
Note that $\frac{d}{dt} \bx_t = f_t \bu_t$.

Then, the chain rule gives the following
\begin{equation} \label{eq:derivative_one}
\begin{aligned}
\frac{d}{dt} \|\bx_t\|^{2m} &= m \|\bx_t\|^{2(m-1)} \frac{d}{dt} \langle \bx_t,\bx_t \rangle \\
& = 2m \|\bx_t\|^{2(m-1)} \Re \langle \bx_t, \frac{d}{dt}\bx_t \rangle \\
& = 2m \|\bx_t\|^{2(m-1)} \Re \langle \bx_t, f_t \bu_t  \rangle
\end{aligned}
\end{equation}
and 
\begin{equation} \label{eq:derivative_two}
\begin{aligned}
&\frac{d}{dt} \left( \|\bx_t\|^{2(m-1)} \bx_t \right) \\
& = \|\bx_t\|^{2(m-1)} \frac{d}{dt} \bx_t + (m-1) \|\bx_t\|^{2(m-2)}\frac{d}{dt} \|\bx_t\|^2 \bx_t \\
& = \|\bx_t\|^{2(m-1)}  f_t \bu_t +2 (m-1) \|\bx_t\|^{2(m-2)} \Re \langle \bx_t, f_t \bu_t \rangle \bx_t
\end{aligned}
\end{equation}

Integrating \eqref{eq:derivative_one}, we obtain
\begin{equation} \label{eq:xT2m}
\begin{aligned}
\|\bx_T\|^{2m} &= \int_0^T 2m \|\bx_t\|^{2(m-1)} \Re \langle \bx_t, f_t \bu_t   \rangle dt \\
& = \int_0^T 2m  \Re \left \langle \int_0^t \frac{d}{d s} \left( \|\bx_s\|^{2(m-1)} \bx_s \right) ds, f_t \bu_t \right \rangle dt    \\
& = 2m \int_0^T   \Re \left \langle \int_0^t   \|\bx_s\|^{2(m-1)}  f_s \bu_s d s, f_t \bu_t \right \rangle dt \\
& \quad + 2m \int_0^T \Re \left\langle \int_0^t 2 (m-1) \|\bx_s\|^{2(m-2)}
\right. \\ & \hspace{50pt} \left.  \Re \langle \bx_s, f_s \bu_s \rangle \bx_s ds, f_t \bu_t \right\rangle dt.
\end{aligned}
\end{equation}
Since $\|\bx_T\|^{2m}$ is a non-negative scalar, $g_T = |g_T|$. Therefore, we take expectation and absolute value on the right of the last \eqref{eq:xT2m} and examine the upper bounds of the two terms separately.

For the first term of \eqref{eq:xT2m}, we have
\begin{equation} \label{eq:first_term}
\begin{aligned}
  &\left| \bbE\left[ \int_0^T   \Re \left \langle \int_0^t   \|\bx_s\|^{2(m-1)}  f_s \bu_s d s, f_t \bu_t \right \rangle dt \right]\right| \\
  &= \left|  \int_0^T  \int_0^t \Re \left( \bbE\left[ \left \langle f_t \bu_t ,  \|\bx_s\|^{2(m-1)}  f_s \bu_s \right \rangle  \right] \right)d s dt \right| \\
  & \le \int_0^T  \int_0^t \left| \bbE\left[ \left \langle f_t \bu_t ,  \|\bx_s\|^{2(m-1)}  f_s \bu_s \right \rangle  \right] \right| d s dt \\
  & =  \int_0^T  \int_0^t |\bar{f_t} f_s| \left| \bbE\left[ \left \langle \bu_t ,  \|\bx_s\|^{2(m-1)}  \bu_s \right \rangle  \right] \right| d s dt.
\end{aligned}
\end{equation}
where the first equality uses the fact that linear operators commute and Fubini's theorem and the first inequality uses triangle inequality and modulus inequality for the real part. 

Let $\bz_s = \|\bx_s\|^{2(m-1)} \bu_s$. By definition, $\bz_s$ is $\cF_s$-measurable. 

Now, we focus on bounding $\left| \bbE\left[ \left \langle \bu_t ,  \|\bx_s\|^{2(m-1)}  \bu_s \right \rangle  \right] \right|$. Firstly, applying Lemma~\ref{lem:innerProd} with $p=2m, q= \frac{2m}{2m-1}$ gives 
\begin{align} \label{eq:expacted_innerprod}
\left| \bbE\left[ \left \langle \bu_t ,  \|\bx_s\|^{2(m-1)}  \bu_s \right \rangle  \right] \right| \le 2 \gamma_{2m}(t-s, \bu) \|\bz_s \|_{L_q}.
\end{align}

Then, applying H\"{o}lder's inequality to $\|\bz_s \|_{L_q} $ gives
\begin{align} \label{eq:z_Lq}
& \bbE\left[ \left\| \|\bx_s\|^{2(m-1)} \bu_s \right\|^q\right]^{1/q} \nonumber \\
& = \bbE\left[  \|\bx_s\|^{2(m-1)q} \|\bu_s\|^q \right]^{1/q} \nonumber\\
& \le \left( \left\| \|\bx_s\|^{2(m-1)q} \right\|_{L_\alpha} \left\| \|\bu_s\|^{q}  \right\|_{L_\beta}\right)^{1/q}  \nonumber \\
& = \bbE\left[ \|\bx_s\|^{2(m-1)q \alpha}  \right]^{\frac{1}{q \alpha}} \bbE\left[ \|\bu_s\|^{q\beta}\right]^{\frac{1}{q \beta}}.
\end{align}
Here, choose $\beta$ such that $q \beta = 2m$, i.e. $\beta = 2m-1$, and so $q \alpha = \frac{2m}{2m-1} \frac{2m-1}{2m-2} = \frac{2m}{2m-2}$.

Therefore, \eqref{eq:z_Lq} is upper bounded by 
\begin{align} \label{eq:z_Lq_bound}
 g_s^{1- \frac{1}{m}} M_{2m}(\bu).
\end{align}

Therefore, plugging \eqref{eq:z_Lq_bound} into \eqref{eq:expacted_innerprod}, we obtain 
\begin{align} \label{eq:first_bound}
\left| \bbE\left[ \left \langle \bu_t ,  \|\bx_s\|^{2(m-1)}  \bu_s \right \rangle  \right] \right|  \le 2 \gamma_{2m}(t-s, \bu) M_{2m}(\bu) g_s^{1- \frac{1}{m}}.
\end{align}

Now, we focus the second term of \eqref{eq:xT2m} and with similar steps as in \eqref{eq:first_term}, we have
\begin{align*}
& \left| \bbE\left[ \int_0^T \Re \left\langle \int_0^t \|\bx_s\|^{2(m-2)}
 \Re \langle \bx_s, f_s \bu_s \rangle \bx_s ds, f_t \bu_t \right\rangle dt \right]\right| \\
 & \le  \int_0^T \int_0^t \left|\bbE\left[ \left\langle f_t \bu_t, \|\bx_s\|^{2(m-2)}
 \Re \langle \bx_s, f_s \bu_s \rangle \bx_s  \right\rangle\right] \right| ds dt.
\end{align*}
Let $\tilde{\bz}_s =  \|\bx_s\|^{2(m-2)}
 \Re \langle \bx_s, f_s \bu_s \rangle \bx_s $ and $\tilde{\bz}_s$ is $\cF_s$-measurable. Recall that $q= \frac{2m}{2m-1}$ and we have
 \begin{align*}
&\bbE\left[\left\| \|\bx_s\|^{2(m-2)}
 \Re \langle \bx_s, f_s \bu_s \rangle \bx_s \right\|^{q} \right]^{1/q} \\
& = \bbE\left[\|\bx_s\|^{2(m-2)q}
\|\Re \langle \bx_s, f_s \bu_s \rangle\|^q \left\|  \bx_s \right\|^{q} \right]^{1/q} \\
& \le \bbE\left[\|\bx_s\|^{2(m-2)q}
\| \bx_s\|^q \|f_s \bu_s \|^q \left\|  \bx_s \right\|^{q} \right]^{1/q} \\
& =  \bbE\left[\|\bx_s\|^{2(m-1)q}
\|f_s \bu_s \|^q \right]^{1/q} \\
& \le \bbE[ \|\bx_s\|^{2(m-1)q \alpha}]^{\frac{1}{q \alpha}} \bbE[
\|f_s \bu_s \|^{q \beta}]^{\frac{1}{q\beta}}
 \end{align*}
where the first equality uses absolute homogeneity, the first inequality uses the Cauchy–Schwarz inequality and the second inequality uses  H\"{o}lder's inequality by choosing $\alpha$, $\beta$ same as before:  $\alpha = \frac{2m-1}{2m-2}$ and $\beta = 2m-1$. 

Note that 
\begin{align*}
\bbE\left[ \| f_s \bu_s \|^{q \beta}\right]^{\frac{1}{q\beta}} \le |f_s| \bbE\left[ \|\bu_s \|^{q \beta}\right]^{\frac{1}{q\beta}}
\end{align*}
and 
\begin{align*}
\gamma_{2m}(t-s, f_t \bu) = |f_t|\gamma_{2m}(t-s, \bu).
\end{align*}

Similarly, using Lemma~\ref{lem:innerProd} and the bounds above, we obtain:
\begin{equation}\label{eq:second_bound}
\begin{aligned} 
& \left|\bbE\left[ \left\langle f_t \bu_t, \|\bx_s\|^{2(m-2)}
 \Re \langle \bx_s, f_s \bu_s \rangle \bx_s  \right\rangle\right] \right| \\
 & \le  2|f_t||f_s| \gamma_{2m}(t-s, \bu) M_{2m}(\bu) g_s^{1-\frac{1}{m}}.
\end{aligned}
\end{equation}
Therefore, combining the bounds in \eqref{eq:first_bound} and \eqref{eq:second_bound}, we have
\begin{align}
  \nonumber
& g_T \le 4m(2m-1) M_{2m}(\bu) \\ \nonumber
& \hspace{30pt} \int_0^T \int_0^t 
|f_s| |f_t| \gamma_{2m}(t-s, \bu) g_s^{1- \frac{1}{m}} ds dt \\ \nonumber
& = 4m(2m-1) M_{2m}(\bu) \\
& \hspace{30pt} \int_0^T g_s^{1- \frac{1}{m}} |f_s|\int_s^T |f_t| \gamma_{2m}(t-s, \bu)  dt ds
\label{eq:g_T_bound_intermediate}
\end{align}
where the equality is obtained via changing the order of the integration. Set $C= 4m(2m-1) M_{2m}(\bu)$.

Define the function $h_s$ as below and then apply the Cauchy–Schwarz inequality
\begin{align*}
h_s &= \int_s^T |f_t| \gamma_{2m}(t-s, \bu)  dt \\
& \le \left( \int_s^T \left( |f_t| \gamma_{2m}(t-s, \bu)^{1/2} \right)^2 dt \right)^{1/2} \cdot \\
& \hspace{60pt} \left( \int_s^T \left( \gamma_{2m}(t-s, \bu)^{1/2} \right)^2 dt \right)^{1/2} \\
& \le \Gamma_{2m}(\bu)^{1/2} \left( \int_s^T \left( |f_t| \gamma_{2m}(t-s, \bu)^{1/2} \right)^2 dt \right)^{1/2} .
\end{align*}

Then, 
\begin{align} 
  \nonumber
& \|h_s\|_{L_2}  \\
\nonumber
&\le  \Gamma_{2m}(\bu)^{1/2}  \left(\int_0^T  \int_s^T  |f_t|^2 \gamma_{2m}(t-s, \bu) dt ds \right)^{1/2}\\
\nonumber
&\le  \Gamma_{2m}(\bu)^{1/2}  \left(\int_0^T  |f_t|^2 \int_0^t   \gamma_{2m}(t-s, \bu) ds dt \right)^{1/2}\\
& \le \Gamma_{2m}(\bu) \left(\int_0^T  |f_t|^2 dt \right)^{1/2} . \label{eq:hs_norm}
\end{align}

Also, applying Lemma~\ref{lem:support_innerProd} to \eqref{eq:g_T_bound_intermediate} by setting $k_s = C |f_s||h_s|$,  
we obtain
\begin{align} \label{eq:g_T_bound}
g_T \le (C/m)^m \left( \int_0^T |f_s| | h_s|ds\right)^m.
\end{align}

Then, by H\"{o}lder's inequality, we have 
\begin{align} \label{eq:int_fs_hs}
\int_0^T |f_s| |h_s|ds \le \left( \int_0^T |f_s|^2 ds \right)^{1/2} \left( \int_0^T |h_s|^2 ds \right)^{1/2}.
\end{align}

Plugging \eqref{eq:hs_norm} into \eqref{eq:int_fs_hs} and further into \eqref{eq:g_T_bound} and then raising both sides of \eqref{eq:g_T_bound} to an exponent $\frac{1}{2m}$ completes the proof.
\QED

The convergence analysis of this paper repeatedly uses the following result, which was given in Corollary 2 in our prior work \cite{zheng2025non}. And it is a discrete-time, multi-dimensional variation of Theorem 1.1 in \cite{gerencser1989class}. The following Corollary uses Lemma~\ref{lem:integral} with $\by_t = \by_{\floor{t}}$, $w_t = w_{\floor{t}}$, and $\Gamma_{2m}(\bu)$ replaced by discrete-time version.  
\begin{corollary}
  {\it
    \label{thm:L_mixing_sum}
    Let $\by_k\in\cY$ be a zero-mean $L$-mixing discrete-time process and let $w_k\in\bbC$. For all $M\ge 1$ and all $q\ge 1$:
    \begin{align*}
    \left\|\sum_{k=0}^{M-1}w_k \by_k \right\|_{L_{2q}} \hspace{-15pt} \le 
      2 \left(2(2q-1)M_{2q}(\by)\Gamma_{d,2q}(\by)
        \sum_{k=0}^{M-1}|w_k|^2
     \right)^{\frac{1}{2}}\hspace{-5pt}.
    \end{align*}
  }
\end{corollary}

\section{Proof of Supporting Lemmas} \label{sec:proofs_supporting}
\subsection{Proof of Lemma~\ref{lem:L-mixing_constant_collection}}
Knowing that $\by[k] - \mu$ is zero-mean and then applying Corollary~\ref{thm:L_mixing_sum} gives 
\begin{align*}
 & \|\hat{\by}_i - h \mu  \|_{L_{2q}} \\
&\le \left\| \sum_{k=0}^{M-1} w_k(s)(\by[i] - \mu) \right\|_{L_{2q}} \\
&\le 2\sqrt{2(2q-1) M_{2q}(\by - \mu) \Gamma_{d,2q}(\by - \mu)}  \sqrt{\sum_{k=0}^{M-1} w_k(s)^2 }.
\end{align*}

We can show that for all $q \ge 1$,
\begin{align} \label{eq:Mbound_shifted}
\|\by - \mu\|_{L_{2q}} \le \|\by  \|_{L_{2q}} +\| \mu \|_{L_{2q}} \le 2M_{2q}(\by). 
\end{align}

This implies that $M_{2q}(\by - \mu) \le 2 M_{2q}(\by)$. Furthermore, directly from the definition, we have 
$
\gamma_{2q}(\cdot, \by -  \mu) =  \gamma_{2q}(\cdot, \by), 
$
which implies that $\Gamma_{d,2q}(\by -  \mu)  = \Gamma_{d,2q}(\by ) $. Then, the desired bound is obtained. 

Bounding $\Gamma_{d,2q}(\tilde{\by})$ is similar to that of Lemma 3 in \cite{zheng2025non}\footnote{The bound on $\Gamma_{d,q}(\hat{\by})$ in Lemma 3 of \cite{zheng2025non} only holds when $\frac{M}{K}=1$ or $2$. The correct general bound in \cite{zheng2025non} should be $2 ( \floor{\frac{M-1}{K}} +1) M_{q}(\by) + (\floor{\frac{M-1}{K}}+1) \Gamma_{d,q}(\by)$. We give the explicit argument here.}. 
%
By construction, $\tilde{\by}_i$ is $\cG_i$ measurable for all $i$ and $\cG_i$ and $\cG_i^+$ are independent. 
  When $\ell K \ge (M-1)$, we have that $(i-\ell)K+(M-1)=\left(iK+k\right)-\left(\ell K+k-(M-1)\right)$, where $\ell K+k-(M-1)\ge 0$ for all $k=0,\ldots,M-1$. In this case, using the triangle inequality, and that $|w_k|\le 1$ gives:
\begin{equation} \label{eq:Gamma_tilde_y_term_2}
  \begin{aligned}[b]
    \MoveEqLeft[0]
    \left\|\tilde{\by}_i-\bbE[\tilde{\by}_i|\cG_{i-\ell}^+] \right\|_{L_{2q}}\\
    &\le 
      \sum_{k=0}^{M-1} \left\|
      \by[iK+k] -\bbE[\by[iK+k]  |\cF_{(i-\ell)K+(M-1)}^+] \right\|_{L_{2q}} \\
    &\le \sum_{k=0}^{M-1}\gamma_{2q}(\ell K +k-(M-1),\by).
  \end{aligned}
\end{equation}  

Now, we can bound $\Gamma_{d,2q}(\tilde{\by})$ via
\begin{align} \label{eq:Gamma_tilde_y}
\MoveEqLeft[0] 
\sum_{\ell=0}^{\infty} \gamma_{2q}(\ell,\tilde{\by}) \le \sum_{l=0}^{\floor{\frac{M-1}{K}}} \gamma_{2q}(\ell,\tilde{\by}) + \sum_{\ell=\floor{\frac{M-1}{K}} + 1}^{\infty} \gamma_{2q}(\ell,\tilde{\by})
\end{align}
We use that 
  $
  \|\tilde{\by}_i-\bbE[\tilde{\by}_i|\cG_{i-\ell}^+]\|_{L_{2q}} \le 2 M_{2q}(\tilde{\by})
  $ for all $0\le \ell \le \floor{\frac{M-1}{K}}$ to bound the first term in \eqref{eq:Gamma_tilde_y}. For the second term in \eqref{eq:Gamma_tilde_y}, we need to plug in the bound from \eqref{eq:Gamma_tilde_y_term_2} and count the repetition of $\gamma_{2q}(\cdot, \by)$ in the summation. It can be shown that one summand, i.e. $\gamma_{2q}(\cdot, \by)$, shows up at most $\floor{\frac{M-1}{K}} + 1$ times in the upper bound of $\gamma_{2q}(\ell, \tilde{\by})$.
(A similar argument for counting the repetitions is given in  the proof of Lemma~\ref{lem:mean_estimate_error}, below.) Therefore, we have the overall bound:
\begin{align*}
\MoveEqLeft[0]
\Gamma_{d,2q}(\tilde{\by}) \\
&\le 2 ( \floor{\frac{M-1}{K}} +1) M_{2q}(\tilde \by) + (\floor{\frac{M-1}{K}}+1) \Gamma_{d,2q}(\by).
\end{align*}
  
Furthermore, bounding $M_{2q}(\tilde{\by} \tilde{\by}^\star) $ is a direct result of Proposition 3 in \cite{zheng2025non} by replacing $\hat{\by}$ by $\tilde{\by}$ there. This gives
\begin{align*}
M_{2q}(\tilde{\by} \tilde{\by}^\star)  &\le M_{4q}(\tilde{\by})^2 \\
\Gamma_{d,2q}(\tilde{\by} \tilde{\by}^\star) &\le 6 M_{4q}(\tilde{\by}) \Gamma_{d,4q}(\tilde{\by}).
\end{align*}

Then, plugging the bounds of $M_{2q}(\tilde{\by})$ and $\Gamma_{d,2q}(\tilde{\by})$ completes the proof.
\hfill \QED

\subsection{Proof of Lemma~\ref{lem:mean_estimate_error}}
The estimate of mean, $\hat{\bmu}_k$, is updated via running average, i.e. $\hat{\bmu}_k = \frac{1}{k} \sum_{i=0}^{k-1} \bar{\by}_i$ for all $k \ge 1$. Recall that $\bar{\by}_i = \frac{1}{M} \sum_{j=0}^{M-1}\by[iK +j]$ and so the total amount of data that $\hat{\bmu}_k$ takes average over is $N = (k-1)K + M$. We can rewrite the running average as $\hat{\bmu}_k = \frac{1}{kM} \sum_{l=0}^{N-1}p[l]\by[l]$, where $p[l]$ counts the number of times that $\by[l]$ appears in the summation.

Note that the last data in $\bar{\by}_i$ is $\by[iK + M-1]$ and the first data in $\bar{\by}_{i+a}$ is $\by[(i+a)K]$. If $\by[\ell]$ appears in the sums for both $\bar{\by}_i$ and $\bar{\by}_{i+a}$, for $a>0$, then
 we must have $iK+M-1 \ge \ell \ge (i+a) K$, i.e. $a \le \frac{M-1}{K}$. Therefore, $\by[l]$ shows up in the summands of $\hat{\bmu}_k$ at most $\floor{a} +1$ times, i.e. $p[l] \le \floor{\frac{M-1}{K}} + 1$.

Therefore, we have the following:
\begin{align*}
\hat{\bmu}_{k} - \mu &= \frac{1}{k} \sum_{i=0}^{k-1} \left(\bar{\by}[i] -  \mu \right)\\
& = \frac{1}{Mk} \sum_{l=0}^{N-1} p[l] \left( \by[l] - \mu \right).
\end{align*}
Taking the $L_{2q}$-norm and applying Theorem~\ref{thm:L_mixing_sum} gives
\begin{align*}
& \left\|\hat{\bmu}_{k} - \mu \right\|_{L_{2q}}\\
 &= \left\|\frac{1}{Mk} \sum_{l=0}^{N-1} p[l]  \left( \by[l] - \mu \right) \right\|_{L_{2q}} \\
& \le 2 \sqrt{2(2q-1) M_{2q}(\by - \mu) \Gamma_{d,2q}(\by - \mu)} \cdot \\
& \hspace{120pt}
\sqrt{\sum_{i=0}^{N -1}\frac{ \left( \floor{\frac{M-1}{K}} + 1 \right)^2}{(Mk)^2}} \\
& \le 2 \sqrt{2(2q-1) M_{2q}(\by - \mu) \Gamma_{d,2q}(\by - \mu)} \frac{\floor{\frac{M-1}{K}} + 1 }{\sqrt{M k}}.
\end{align*}
where the last inequality uses: $K \le M \Rightarrow N \le kM$.

Then applying the bound from \eqref{eq:Mbound_shifted} and the definition of $\Gamma_{d,q}(\by)$ completes the proof.
%
%
\hfill \QED

\section{Proofs of Main Results} \label{sec:proofs_main}

\subsection{Proof of Theorem~\ref{theorem_main_batch}}
From the construction of the batch algorithm \eqref{eq:alg_batch}, we have
\begin{equation} \label{eq:decomposition_batch}
\begin{aligned}[b] 
\hat{\bPhi}_k - \bar{\Phi} &= \frac{1}{k} \sum_{i=0}^{k-1} (\hat{\by}_i - h \hat{\bmu}_k)(\hat{\by}_i - h \hat{\bmu}_k)^\star -\bar{\Phi} \\
& =\frac{1}{k} \sum_{i=0}^{k-1}  \left( (\hat{\by}_i - h \mu) (\hat{\by}_i - h \mu)^\star  -\bar{\Phi} \right)  \\ 
  & \hspace{10pt} 
 +  \frac{1}{k} \sum_{i=0}^{k-1} (\hat{\by}_i - h \mu) ( h \mu - h \hat{\bmu}_k)^\star  \\
 &  \hspace{10pt}   + \frac{1}{k} \sum_{i=0}^{k-1} ( h \mu - h \hat{\bmu}_k) (\hat{\by}_i - h \mu)^\star   \\
 & \hspace{10pt}  + \frac{1}{k} \sum_{i=0}^{k-1} ( h \mu - h \hat{\bmu}_k) ( h \mu - h \hat{\bmu}_k)^\star  
\end{aligned} 
\end{equation}
where the second equality comes from adding and subtracting $h \mu$ from $\hat{\by}_i - h \hat{\bmu}_k$ and then expanding the terms.

Now, we bound the $L_{2q}$-norm of the four terms in \eqref{eq:decomposition_batch} separately.

For the first term, Corollary~\ref{thm:L_mixing_sum} gives
\begin{equation} \label{eq:bound_first_term_batch}
\begin{aligned}[b]
& \left\| \frac{1}{k} \sum_{i=0}^{k-1}  \left( (\hat{\by}_i - h \mu) (\hat{\by}_i - h \mu)^\star  -\bar{\Phi} \right) \right\|_{L_{2q}} \\
 &  \le  2 \sqrt{2(2q-1)M_{2q}(\tilde{\by}\tilde{\by}^* - \bar{\Phi}) \Gamma_{d,2q}(\tilde{\by}\tilde{\by}^\star- \bar{\Phi})} \sqrt{ \sum_{i=0}^{k-1} \frac{1}{k^2}} \\
 & \le 4 \sqrt{(2q-1)M_{2q}(\tilde{\by}\tilde{\by}^\star) \Gamma_{d,2q}(\tilde{\by}\tilde{\by}^\star)} \frac{1}{\sqrt{k}}.
\end{aligned}
\end{equation}

Note that similar to \eqref{eq:Mbound_shifted}, applying the triangle inequality and then Jensen's inequality gives $M_{2q}(\tilde{\by}\tilde{\by}^\star - \bar{\Phi}) \le 2 M_{2q}(\tilde{\by}\tilde{\by}^\star) $. Besides, $\Gamma_{d,2q}(\tilde{\by}\tilde{\by}^\star- \bar{\Phi}) = \Gamma_{d,2q}(\tilde{\by}\tilde{\by}^\star)$ holds directly from the definition of the $L$-mixing property.

The second and third terms can be bounded in the same way:
\begin{equation} \label{eq:bound_cross_batch} 
\begin{aligned}[b]
&\left\| \frac{1}{k} \sum_{i=0}^{k-1} (\hat{\by}_i - h \mu) ( h \mu - h \hat{\bmu}_k)^\star \right\|_{L_{2q}} \\
&\le \frac{1}{k} \sum_{i=0}^{k-1} \|\hat{\by}_i - h \mu\|_{L_{4q}} \|h \mu - h \hat{\bmu}_k\|_{L_{4q}} \\
&\le \frac{1}{k} M_{4q}(\tilde{\by}) \sqrt{M} c_{2q} \frac{ 1}{\sqrt{M}} \sum_{i=0}^{k-1} \frac{1}{\sqrt{k}} \\
& \le M_{4q}(\tilde{\by})c_{2q}  \frac{1}{\sqrt{k}}
\end{aligned}
\end{equation}
where the first inequality uses the absolute homogeneity of $L_q$-norm, the triangle inequality followed by the fact that for any vectors $\|xy^\star\|_F = \|x\|_2 \|y\|_2$ and the Cauchy-Schwarz inequality.
The third inequality uses Lemma~\ref{lem:mean_estimate_error}.
Note that $|h| = \sqrt{M} $ from the Cauchy-Schwarz inequality.

Similarly, the fourth term is bounded as:
\begin{equation} \label{eq:bound_quadratic_batch}
\begin{aligned}[b]
& \left\|\frac{1}{k} \sum_{i=0}^{k-1} ( h \mu - h \hat{\bmu}_k) ( h \mu - h \hat{\bmu}_k)^\star  \right\|_{L_{2q}} \\
& \le \frac{1}{k} \sum_{i=0}^{k-1} \|h \mu - h \hat{\bmu}_k\|^2_{L_{4q}} \\
& \le \frac{1}{k} M c_{2q}^2 \frac{1}{M} \sum_{i=0}^{k-1} \frac{1}{k}  \le c_{2q}^2  \frac{1}{k}.
\end{aligned}
\end{equation}

Therefore, combining the bounds in \eqref{eq:bound_first_term_batch}, \eqref{eq:bound_cross_batch}, and \eqref{eq:bound_quadratic_batch} gives 
\begin{align*}
\left\|\hat{\bPhi}_k - \bar{\Phi}  \right\|_{L_{2q}} &\le \left(  4 \sqrt{(2q-1)M_{2q}(\tilde{\by}\tilde{\by}^\star) \Gamma_{d,2q}(\tilde{\by}\tilde{\by}^\star)} \right. \\
& \hspace{10pt} \left. +  2 M_{4q}(\tilde{\by})  c_{2q}  \right) \frac{1}{\sqrt{k}} + c_{2q}^2 \frac{1}{k}.
\end{align*}

Then, plugging the constant bounds from Lemma~\ref{lem:L-mixing_constant_collection} followed by the monotonicity of $L_q$-norm completes the proof.
\hfill \QED

\subsection{Proof of Theorem~\ref{theorem_main_online}}
The proof largely remains the same as that of Theorem~\ref{theorem_main_batch} and only varies in the parts containing the time-varying mean estimate $\hat{\bmu}_i$.


From the construction of the online algorithm \eqref{eq:alg_online}, we have
\begin{align}
\nonumber
\hat{\bPhi}_k - \bar{\Phi} &= \frac{1}{k} \sum_{i=0}^{k-1} ( \Delta \hat \by_i \Delta \hat \by_i^* - \bar{\Phi}) \\
\nonumber
& =\frac{1}{k} \sum_{i=0}^{k-1}  \left( (\hat{\by}_i - h \mu) (\hat{\by}_i - h \mu)^*  -\bar{\Phi} \right)  \\ 
\nonumber
  & \hspace{10pt} 
 +  \frac{1}{k} \sum_{i=0}^{k-1} (\hat{\by}_i - h \mu) ( h \mu - h \hat{\bmu}_i)^*  \\
\nonumber
 &  \hspace{10pt}   + \frac{1}{k} \sum_{i=0}^{k-1} ( h \mu - h \hat{\bmu}_i) (\hat{\by}_i - h \mu)^*   \\
 & \hspace{10pt}  + \frac{1}{k} \sum_{i=0}^{k-1} ( h \mu - h \hat{\bmu}_i) ( h \mu - h \hat{\bmu}_i)^*  \label{eq:decomposition}
\end{align}
where the second equality comes from plugging 
$\Delta \hat \by_k = \hat{\by}_k - h\mu+ h \mu - h \hat{\bmu}_k.$

Similarly, we bound the $L_{2q}$-norm of the four terms in \eqref{eq:decomposition} separately. The bound of the first term in \eqref{eq:decomposition} remains exactly the same. The processes of getting the bounds on the remaining three terms are slightly different.

The second and third term share the same bound:
\begin{equation} \label{eq:bound_cross_online}
\begin{aligned}[b]
&\left\| \frac{1}{k} \sum_{i=0}^{k-1} (\hat{\by}_i - h \mu) ( h \mu - h \hat{\bmu}_i)^\star \right\|_{L_{2q}} \\
& \le \frac{1}{k} \sum_{i=1}^{k-1} \|\hat{\by}_i - h \mu\|_{L_{4q}} \|h \mu - h \hat{\bmu}_i\|_{L_{4q}} \\
& \hspace{10pt} + \frac{1}{k} \|\hat{\by}_0 - h \bar{\mu}\|_{L_{4q}} \|h \mu - h \hat{\bmu}_0\|_{L_{4q}}
 \\
&\le \frac{1}{k} M_{4q}(\tilde{\by}) c_{2q} \sum_{i=1}^{k-1} \frac{1}{\sqrt{i}} + \frac{1}{k} M_{4q}(\tilde{\by}) \sqrt{M} \|\mu\|_{L_{4q}} \\
& \le 3 M_{4q}(\tilde{\by}) c_{2q}\frac{1}{\sqrt{k}} + M_{4q}(\tilde{\by}) M_{4q}(\by)  \frac{\sqrt{M} }{k}
\end{aligned}
\end{equation}
where the last two inequalities use $\hat{\bmu}_0 = 0$ and $\|\mu\|_{L_{4q}} = M_{4q}(\by)$ and the Riemann sum bound:  $\sum_{i=1}^{k-1} \frac{1}{\sqrt{i}} \le 1 + \int_{1}^{k-1}\frac{1}{\sqrt{t}}dt \le  1+ 2\sqrt{k} \le 3 \sqrt{k} $ for all $k\ge 1$.  

Similarly, the fourth term is bounded as below:
\begin{align} \label{eq:bound_quadratic_online}
\nonumber& \left\|\frac{1}{k} \sum_{i=0}^{k-1} ( h \mu - h \hat{\bmu}_i) ( h \bar{\mu} - h \hat{\bmu}_i)^*  \right\|_{L_{2q}} \\
\nonumber & \le c_{2q}^2 \frac{1}{k} \sum_{i=1}^{k-1}  \frac{1}{i} + \frac{1}{k} M \|\mu\|_{L_{4q}}^2 \\
& \le  2  c_{2q}^2 \frac{1}{\sqrt{k}} +  M_{4q}(\by)^2 \frac{M}{k}
\end{align}
where the last inequality uses the following approximations: $\sum_{i=1}^{k-1}  \frac{1}{i} = 1 + \int_1^{k-1} \frac{1}{t} dt \le 1+ \log (k-1)$, $\log (k -1) \le \log k \le \sqrt{k}$  for all $k \ge 2$ and $1 \le \sqrt{k}$ for all $k \ge 1$.

Therefore, combining the bounds in \eqref{eq:bound_first_term_batch}, \eqref{eq:bound_cross_online}, and \eqref{eq:bound_quadratic_online} gives
\begin{align*}
&\left\|\hat{\bPhi}_k - \bar{\Phi}  \right\|_{L_{2q}} \\
& \le \left( 4 \sqrt{(2q-1) M_{2q}(\tilde{\by}\tilde{\by}^* ) \Gamma_{d,2q}(\tilde{\by}\tilde{\by}^*)}  + 6  M_{4q}(\tilde{\by})c_{2q}  \right. \\
& \hspace{10pt}\left.   
 +  2 c_{2q}^2 \right) \frac{1}{\sqrt{k}} + \left(2 M_{4q}(\tilde{\by})M_{4q}(\by) +   M_{4q}(\by)^2\right) \frac{M}{k}.
\end{align*}



Then, plugging the constant bounds from Lemma~\ref{lem:L-mixing_constant_collection} followed by the monotonicity of $L_q$-norm completes the proof.
\hfill \QED
\subsection{Proof of Corollary~\ref{cor:known_mean}}

When $\mu$ is known, the last three terms of \eqref{eq:decomposition_batch} or \eqref{eq:decomposition} vanish because $\hat{\bmu}_k - \mu = 0$ for all $k \ge 1$. Then, plugging the bounds in Lemma~\ref{lem:L-mixing_constant_collection} completes the proof.
\hfill \QED

\section{Simulation} \label{sec:sim}

The simulation setup remains the same as \cite{zheng2025non}. Data samples are generated from measurements of the same finite-state Markov chain. The Markov states are in $\{0,1\}$ and the transition matrix is $
P = \begin{bmatrix}
0.3 & 0.7 \\
0.5 & 0.5
\end{bmatrix}.
$ 
Note the new bound does not require the zero-mean assumption, so we don't have to demean the data as in the previous work \cite{zheng2025non}.

The work \cite{ gerencser2002new} shows that the Markov chain in our example is $L$-mixing and the Doeblin coefficient is $\delta = 0.72$. Therefore, we have the following $L$-mixing statistics: $\Gamma_{d,4q}(\by) \le   4 G_{\max} \frac{1}{1- (1-\delta)^{\frac{1}{4q}}} \le \frac{4 G_{\max}}{\delta} 4q$ and $M_{4q}(\by) \le G_{\max}$. 

\begin{figure}[t]
      \centering
      \subfloat[\textbf{Bartlett Estimator}. $M =5, L = 10^7$]{
      \includegraphics[width=.7\columnwidth]{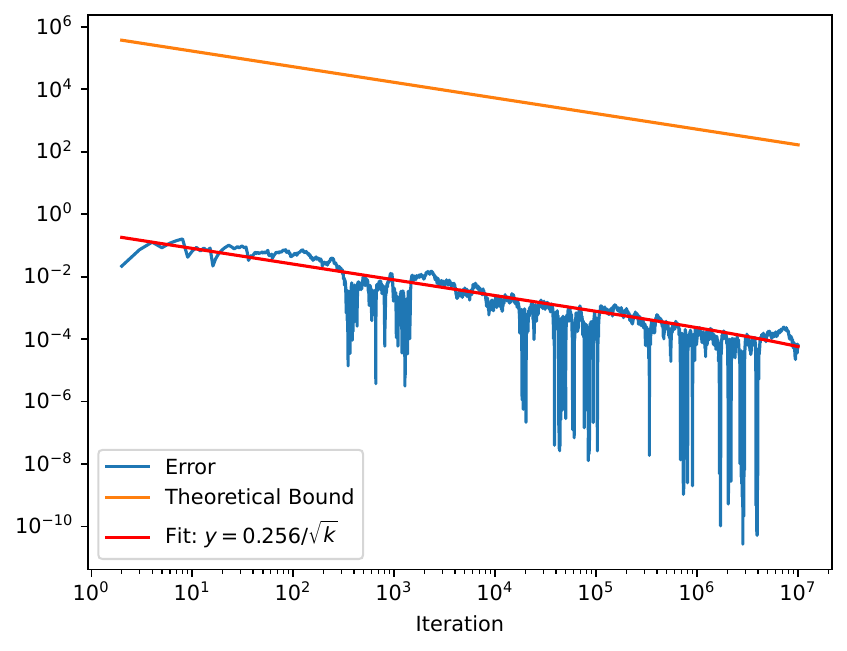}
    
      }
      \hfill
      \subfloat[\textbf{Welch Estimator}. Hann Window, $M =16, K=8, L = 10^7$]{
      \includegraphics[width=.7\columnwidth]{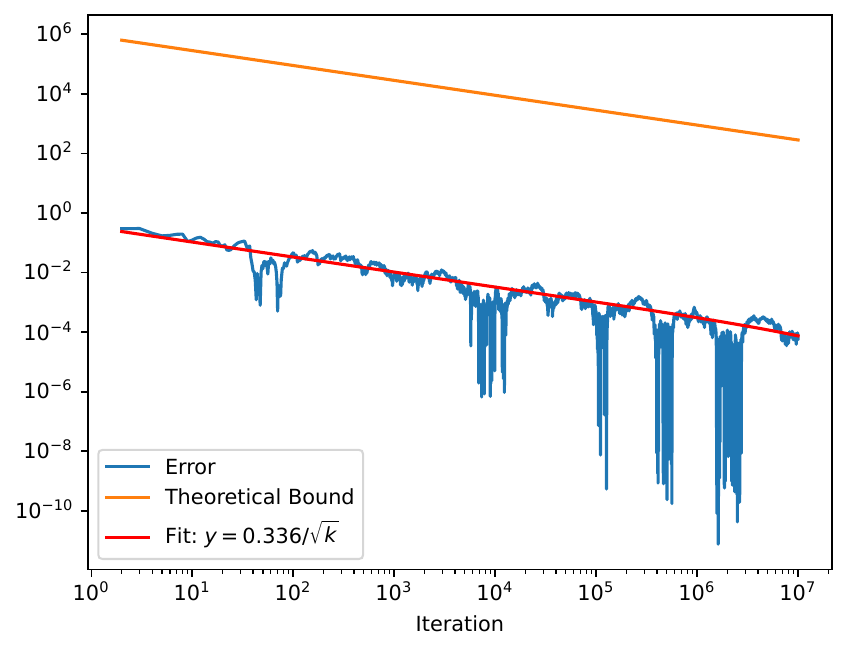}
      
      }
	\caption{Concentration of estimate to its mean on finite Markov chain data}
	\label{fig:errors}
  \end{figure}

We only present the simulation result for online algorithms since batch algorithms should have similar phenomena but requires many runs of simulations. The theoretical bounds in Fig. \ref{fig:errors} is from Theorem~\ref{thm:highProb} with $\nu=0.1$ regarding the online algorithm. Note that the iteration in the figures refer to the number of data segments used in the algorithms, namely $k$. The empirical errors are calculated via $\|\hat \bPhi_k(s) -\bar{\Phi}(s)\|$ and are well below the theoretical bounds. 
Additionally, the best-fit curves of the empirical errors for both estimators show a $\Theta(\frac{1}{\sqrt{k}})$ dependence on $k$. We can observe that the errors evolve almost the same as those in \cite{zheng2025non} other than the first few steps. This is because the mean estimate converges to true average of the data very quickly. It is clear that the theoretical bounds are quite conservative. This is due to large constant factors. Note here, that we are giving the \emph{first} constant factors for this problem, and we are giving explicit bounds on them. (In many cases, bounds are stated in terms of \emph{unquantified} constant factors.)  Tightening these bounds is an important area for future research.

\section{Conclusion and Future Work} \label{sec:conclusion}

This work derives non-asymptotic error bounds for Bartlett and Welch estimators for $L$-mixing data with unknown means using  batch and online algorithms. 
High probability error bounds are also obtained, and we have simulated a finite Markov chain to verify the theory. Our error bounds are $O(\frac{1}{\sqrt{k}})$, where $k$ is the number of data segments used in the algorithm, which are tighter than the results obtained in \cite{lamperski2023nonasymptotic} and \cite{zheng2025non}. 



One future direction is to obtain tighter error bounds, possibly by conducting a frequency-dependent analysis.
Furthermore, different choices of step size may improve the algorithm performance under the stationary assumption, and may also allow analysis in non-stationary settings.
Such work would benefit the analysis of complex dynamics where the presence of nonstationary time-series data is unavoidable. 

\printbibliography
\end{document}